\hfuzz=6pt
\font\titlefont=cmbx10 scaled\magstep1

\magnification=1200
\line{}
\vskip 1.5cm
\centerline{\titlefont A MODEL FOR THE CONTINUOUS}
\smallskip
\centerline{\titlefont q-ULTRASPHERICAL POLYNOMIALS}
\vskip 2cm
\centerline{\bf Roberto Floreanini}
\smallskip
\centerline{Istituto Nazionale di Fisica Nucleare, Sezione di Trieste}
\centerline {Dipartimento di Fisica Teorica,
Universit\`a di Trieste}
\centerline{Strada Costiera 11, 34014 Trieste, Italy}
\vskip 1cm
\centerline{\phantom{$^{(*)}$}{\bf Luc Vinet}\footnote{$^{(*)}$}{Supported 
in part by the National Sciences and Engineering
Research Council \hbox{(NSERC)} of Canada and the Fonds FCAR of Qu\'ebec.}}
\smallskip
\centerline{Centre de Recherches Math\'ematiques}
\centerline{Universit\'e de Montr\'eal}
\centerline{Montr\'eal, Canada H3C 3J7}
\vskip 1.5cm
\centerline{\bf Abstract}
\smallskip\midinsert\narrower\narrower\noindent  
We provide an algebraic interpretation for two classes of continuous
$q$-polynomials. Rogers' continuous $q$-Hermite polynomials and
continuous $q$-ultraspherical polynomials are shown to realize,
respectively, bases for representation spaces of the $q$-Heisenberg
algebra and a $q$-deformation of the Euclidean algebra in these 
dimensions. A generating function for the continuous
$q$-Hermite polynomials and a $q$-analog of the Fourier-Gegenbauer
expansion are naturally obtained from these models.
\endinsert

\vfill\eject

{\bf 1. INTRODUCTION}
\bigskip

The algebraic theory of $q$-special functions$^{1,2}$ is
currently being actively developed.  Indeed, it has been realized
that quantum groups and algebras offer a unifying framework for
describing and studying these functions.  In a series  of papers 
on this topic,$^{3-12}$ we, among others, have shown that many
$q$-orthogonal polynomials and $q$-functions appear as matrix
elements or basis vectors of quantum algebra representations.  We
have further used these observations to derive and obtain in a natural way, 
various relations and identities that $q$-special
functions obey.  In most cases, the $q$-polynomials that we
encountered, were orthogonal with respect to a discrete measure. 
There are, however, numerous sets of $q$-polynomials orthogonal with 
respect to continuous measures.$^{1,2}$  Two classes of these
continuous $q$-orthogonal polynomials will be considered here and
given algebraic interpretations.

Our approach to basic special functions can be illustrated  by
considering the connection between $q$-Bessel functions 
and the two-dimensional quantum Euclidean algebra 
${\cal U}_q\bigl(E(2)\bigr)$.$^{4,5}$ This Hopf algebra is generated by
 the elements $D, D^{-1}$, $P_+$ and $P_-$ satisfying the defining
relations:
$$
[P_+, P_-] =0, \qquad DP_\pm = q^{\pm 1} P_\pm D, \qquad DD^{-1} =
D^{-1}D = 1\ .\eqno(1.1)
$$
The coproduct $\Delta$, antipode $S$ and counit $\varepsilon$ are
determined by
$$
\eqalign{&\Delta(D) = D \otimes D\ ,\cr
         &S(D)= D^{-1}\ ,\cr
         &\varepsilon(D)=\varepsilon(1) = 1\ ,\cr}\qquad
\eqalign{&\Delta(P_\pm) = P_\pm \otimes D^{-1/2} + D^{1/2} \otimes P_\pm\ ,\cr
         &S(P_\pm)=-q^{\mp 1/2} P_\pm\ ,\cr
         &\varepsilon(P_\pm) = 0\ .}\eqno(1.2)
$$
Representations of ${\cal U}_q\bigl(E(2) \bigr)$ are characterized by
two complex numbers $\omega$ and $m_0$, with $\omega \neq 0$ and 
$0\leq {\rm Re}\, m_0 < 1$.  The representation spaces have for basis, vectors
$f_m$, with the index $m$ running over the elements
 of the set $S = \{ m_0 + n \colon n \in{\bf Z}\}$.  In these bases
the generators act as follows:
$$
P_\pm f_m
	= \omega f_{m \pm 1}\ ,\qquad
D^{\pm 1} f_m 
	= q^{\pm m} f_m\ .\eqno(1.3)
$$
Consider now the following two $q$-analogs of the exponential
function:$^1$
$$
\eqalign{
&e_q(z)= \sum_{n=0}^\infty {1\over(q;q)_n} z^n =
{1\over(z;q)}_\infty\ ,\cr
&E_q(z)= \sum_{n=0}^\infty {q^{n(n-1)/2}\over(q;q)_n} z^n =
(-z;q)_\infty\ .}\eqno(1.4)
$$
The symbols $(a;q)_\alpha$ are referred to as $q$-shifted factorials,
and are defined by
$$
(a;q)_\alpha={(a;q)_\infty\over(aq^\alpha;q)_\infty}\ , \eqno(1.5a)
$$
with
$$
(a;q)_\infty=\prod_{k=0}^\infty (1 - aq^k)\ ,  \qquad |q| < 1\ .\eqno(1.5b) 
$$

The $q$-exponentials $e_q(z)$ and $E_q(z)$ are eigenfunctions of the
$q$-derivative operators $D_z^+$ and $D_z^-$ respectively.  These
operators are defined by 
$$
D_z^\pm = {1\over z} (1 - T_z^{\pm 1}), \eqno(1.6)
$$
where $T_z$ is the $q$-shift operator $T_z\, f(z) = f(qz)$.
It is easy to check that
$$
D_z^+ e_q(\lambda z) = \lambda e_q(\lambda z)\ ,\qquad
D_z^- E_q(\lambda z)= -q^{-1} \lambda E_q(\lambda z)\ .\eqno(1.7)
$$
Note also that $\lim_{q \rightarrow 1^-} e_q\bigl( z(1-q) \bigr) = 
\lim_{q\rightarrow 1^-} E_q \bigl( z(1-q) \bigr) = e^z$.  

In analogy with ordinary Lie theory, we introduce for instance, the
following element in the completion of ${\cal U}_q\bigl( E(2) \bigr)$:
$$
U(\alpha, \beta) = E_q(\alpha P_+)\, E_q (\beta P_-)\ . \eqno(1.8)
$$
In the limit $q \rightarrow 1^-$, $U \bigl[ (1-q)\alpha, (1-q)\beta \bigr]$
goes into the element $e^{\alpha P_+ + \beta P_-}$ of the
two-dimensional Euclidean Lie group.  The matrix elements
$U_{kn}(\alpha, \beta)$ of $U(\alpha, \beta)$ are defined through
$$
U(\alpha, \beta)\, f_{m_0+n} = \sum_{k=-\infty}^\infty 
U_{kn}(\alpha,\beta)\, f_{m_0+k}\ . \eqno(1.9)
$$
Given the action (1.3) of the generators, these matrix elements are
straightforwardly computed and one finds$^{4,5}$
$$
U_{kn} (\alpha, \beta) = q^{(k-n)^2/2} \left( -
{\alpha\over\beta}\right)^{(k-n)/2}\ J_{k-n}^{(2)} \left( 2\,\omega
\left( -{\alpha\beta\over q}\right)^{1/2}; q \right)\ .  \eqno(1.10)
$$
The functions $J_\nu^{(2)}(z;q)$ are $q$-analogs of the Bessel
functions $J_\nu(z)$.  They are defined by the series$^{1,13}$
$$
J_\nu^{(2)} (z;q) = \sum_{n=0}^\infty q^{n(n+\nu)}
{(-1)^n\over(q;q)_n\,(q;q)_{n+\nu}} \left({z\over2}\right)^{2n+\nu}\ ,
\eqno(1.11)
$$
and one readily checks that 
$\lim_{q \rightarrow 1^-} J_\nu^{(2)} \bigl(z(1-q); q \bigr) = J_\nu(z)$.  
This relation between the $q$-Bessel
functions and ${\cal U}_q\bigl( E(2) \bigr)$ proves rather fruitful. 
It is for instance possible to construct a two-variable
realization of ${\cal U}_q \bigl( E(2) \bigr)$ where the matrix
elements $U_{kn} (\alpha, \beta)$, appear as basis vectors of the
associated module.  Within this model, (1.9) is shown to
entail a $q$-analog of Graf's addition formula for Bessel functions.$^{5}$  
Many properties of the $q$-Bessel functions can be
obtained and interpreted in this fashion.  Let us record for future
reference the three-term recurrence relation that the functions
$J_\nu^{(2)}(z;q)$ obey$^{1,13}$
$$
q^\nu\, J_{\nu+1}^{(2)} (z;q) = {2\over z} (1-q^\nu)\, J_\nu^{(2)}(z;q)
- J_{\nu-1}^{(2)} (z;q)\ , \eqno(1.12)
$$
and the asymptotic behavior of these functions as $\nu$ goes to
infinity$^{14}$
$$
J_\nu^{(2)} (z;q) \sim {(z/2)^\nu\over (q;q)_\infty}\ .\eqno(1.13)
$$

In summary, our algebraic interpretation of $q$-special functions
proceeds as follows.  Using certain $q$-exponentials of $q$-algebra
generators, $q$-analogs of Lie group elements are formed and their
matrix elements in representation spaces of the deformed 
algebra are shown to involve $q$-special functions.  Various
models and realizations are then constructed and called upon to
exploit these results.  This approach has been applied to many
situations, but so far, has not been used much to discuss the 
interesting continuous $q$-orthogonal polynomials.  Encompassing all of
these, are the Askey-Wilson polynomials $p_n(x;a,b,c,d \vert q)$. 
They form a four-parameter family and are defined by$^{15,1,2}$
$$
\eqalign{p_n(x;a,b,c,d \vert q) = &(ab,ac,ad;q)_n\, a^{-n}\cr
&\times {}_4\phi_3
\bigg({q^{-n},\atop\ }{abcd\, q^{n-1},\atop ab,}{ae^{i\theta},\atop ac,}
{ae^{-i\theta}\atop ad}\bigg| q; q\bigg)\ .}\eqno(1.14)
$$
We are using the notation
$$
(a_1, a_2, \dots , a_k; q)_\alpha = (a_1; q)_\alpha\, 
(a_2; q)_\alpha \dots (a_k; q)_\alpha\ , \eqno(1.15)
$$
$$
\eqalign{{}_r\phi_s \bigg({a_1, a_2, \dots , a_r\over
b_1, \dots b_s}\bigg| q; z\bigg)
= \sum_{n=0}^\infty & {(a_1, \dots a_r; q)_n\over(q,b_1, \dots ,b_s; q)_n}\cr
&\hskip 1cm\times \left[ (-1)^n q^{n(n-1)/2} \right]^{1+s-r}\ z^n\ .} 
\eqno(1.16)
$$
When $a,b,c,d$ are real, the Askey-Wilson polynomials are orthogonal
over the interval $0 < \theta < \pi$ with respect to the continuous
measure
$$
w(\cos \theta; a,b,c,d) =
\left| {(e^{2i\theta};q)_\infty\over
(ae^{i\theta},be^{i\theta}, ce^{i\theta}, de^{i\theta};q)_\infty}\right|^2
\ .\eqno(1.17)
$$
Particular cases of interest are obtained for special choices of the
parameters.  The continuous $q$-ultraspherical polynomials 
$C_n(x;\beta \vert q)$ of Rogers are 
obtained by setting $a = \beta^{1/2}$,
$b = \beta^{1/2} q^{1/2}$, $c = -\beta^{1/2}$, $d = -\beta^{1/2}q^{1/2}$ 
and changing also the normalization:$^{1,2,17}$
$$
\eqalign{C_n(x; \beta \vert q)= &
{(\beta^2; q)_n\over(\beta q^{1/2}, -\beta, -\beta q^{1/2}, q;q)_n}\cr
&\hskip 2cm\times p_n \bigl( x; \beta^{1/2}, \beta^{1/2}q^{1/2}, -\beta^{1/2},
-\beta^{1/2}q^{1/2} \vert q \bigr)\ .}\eqno(1.18)
$$
By putting $a=0$, $b=0$, $c=0$, $d=0$,  one arrives at the continuous
$q$-Hermite polynomials $H_n(x \vert q)$:$^{1,2,17}$
$$
H_n(x \vert q) = p_n (x; 0,0,0,0 \vert q)\ . \eqno(1.19)
$$
In the following, our attention will be focused on these two classes
of polynomials.  Note that$^{2}$
$$\lim_{q \rightarrow 1^-} C_n (x; q^\lambda ;q) = C_n^\lambda (x)\ ,
\eqno(1.20a)$$
$$\lim_{q \rightarrow 1^-} \left({1-q\over2}\right)^{-n/2} 
H_n \bigg(x\sqrt{{1-q\over2}} \bigg| q \bigg) = H_n(x)\ ,\eqno(1.20b)
$$
where $C_n^{(\lambda)}(x)$ and $H_n(x)$ respectively denote the
ultraspherical and Hermite polynomials.$^{18}$  Let us point out
that neither the $q$-exponential $e_q(x)$, nor $E_q(x)$ proved
relevant in the algebraic interpretation of the polynomials 
$C_n(x; \beta \vert q)$ and $H_n(x \vert q)$.  In fact, it is the
$q$-exponential introduced recently in Ref.[14] that turns out to be
the appropriate $q$-analog of the exponential to use here.  The
properties of this function will be discussed in Section 3, 
after the relation between the continuous $q$-Hermite polynomials
and the $q$-oscillator algebra will have been established in
Section~2.

A generating function for these continuous $q$-Hermite polynomials
will be derived in Section~4.  Next, in Section~5, a $q$-deformation
of the three-dimensional Euclidean algebra will be realized on the
complex 2-sphere.  This construct will be seen to provide 
a nice framework for studying the properties of the continuous
$q$-ultraspherical polynomials.  It will  allow us to present in
Section~6, an algebraic derivation of a $q$-analog (see (6.22)) of the
Fourier-Gegenbauer expansion:$^{14}$
$$
e^{i\vec k \cdot \vec r} = \Gamma(\nu)
\left({kr\over2}\right)^{-\nu} \sum_{n=0}^\infty i^n (\nu + n)\,
J_{\nu + n}(kr)\, C_n^\nu(\cos \theta)\ ,\eqno(1.21)
$$
where $\vec k \cdot \vec r = kr \cos \theta$.  Concluding remarks
will then end the paper.

\vskip 2cm

{\bf 2. CONTINUOUS q-HERMITE POLYNOMIALS AND THE}\hfill\break\indent
{\bf \phantom{2. }q-HEISENBERG ALGEBRA}
\bigskip

The continuous $q$-Hermite polynomials $H_n(x \vert q)$ given in
(1.19) as special cases of the Askey-Wilson polynomials, can also be
defined as follows:$^{2}$
$$
H_n(x \vert q) = \sum_{k=0}^n {(q;q)_n\over (q;q)_k\, (q;q)_{n-k}}\,
e^{i(n-2k)\theta}\ , \qquad x = \cos \theta\ . \eqno(2.1)
$$
Let
$$
z = e^{i\theta}\ , \eqno(2.2)
$$
and introduce the divided difference operators$^{19}$
$$
\eqalignno{
&\tau = {1\over z - z^{-1}}\,  \big(T_z^{1/2} - T_z^{-1/2}\big)\ , &(2.3)\cr
&\tau^*= {q^{-1/2}\over z - z^{-1}} \left({1\over z^2}\, T_z^{1/2} -
z^2\, T_z^{-1/2} \right)\ . &(2.4)}
$$
These operators will be taken to act on functions of $x = (z +z^{-1})/2$.  
Notice that as $q \rightarrow 1$:
$$
\lim_{q \rightarrow 1^-} {2\over(q^{1/2} - q^{-1/2})} \, \tau =
{d\over dx}\ , \eqno(2.5)
$$
$$
\lim_{q \rightarrow 1^-} \tau^*/2 = x\ . \eqno(2.6)
$$
Set 
$$
f_n(x) = H_n (x \vert q), \qquad n=0,1,2,\dots\ . \eqno(2.7)
$$
One can directly verify from the definition (2.1), that
$$
\eqalignno{
&\tau\, f_n(x)= q^{n/2} (1 - q^{-n})\, f_{n-1}(x)\ , &(2.8)\cr
&\tau^*\, f_n(x)=-q^{-(n+1)/2}\, f_{n+1}(x)\ . &(2.9)}
$$
It follows that $\tau$ and $\tau^*$ satisfy
$$
\tau^* \tau - q\, \tau \tau^* = -(1-q)\ . \eqno(2.10)
$$
These operators hence provide a realization of the $q$-Heisenberg
algebra with the continuous $q$-Hermite polynomials occurring as
basis vectors for a representation space of this algebra.  The action
on that basis of the operator$^{19}$
$$
\mu ={1\over z - z^{-1}} \left( -{1\over z}\,  T_z^{1/2} + 
z\, T_z^{-1/2} \right)\ , \eqno(2.11)
$$
is also readily determined:
$$
\mu\, f_n(x) = q^{-n/2}\, f_n(x)\ . \eqno(2.12)
$$
Remark that $\mu$ reduces to the identity operator when $q \rightarrow 1$ and
that
$$
\eqalignno{
&\tau \mu - q^{-1/2} \mu \tau=\, 0\ , &(2.13)\cr
&\tau^* \mu - q^{1/2} \mu \tau^* =\, 0\ . &(2.14)}
$$
It is also natural to consider the operator multiplication by $x$. 
This operator is not independent of $\tau$, $\tau^*$ and $\mu$ since
$$
x\, \mu = -\tau - q^{1/2}\, \tau^*\ . \eqno(2.15)
$$
Its action on the basis $\{ f_n \}$ can be gotten from this last
relation and is in fact embodied in the three-term recurrence
relation of the continuous $q$-Hermite polynomials:$^{1,2}$
$$
2x\, f_n(x) = f_{n+1}(x) + (1 - q^n)\, f_{n-1}(x)\ . \eqno(2.16)
$$
We now want to follow the approach that we outlined in the
Introduction, that is, take $q$-exponentials of algebra elements and
compute their action on representation spaces to arrive at special
function identities.  To this end, we need $q$-exponentials 
that behave nicely under the action of the elementary difference
operator of the models under consideration.  In the present case, we
would thus wish for a $q$-exponential that is an eigenfunction of the
divided difference operator $\tau$.  Such a function 
has been introduced recently and will be the object of the next
section.

\vskip 2cm

{\bf 3. THE q-EXPONENTIAL ${\cal E}_q$}
\bigskip

Consider the divided difference operator $\tau$ given in (2.3).  We
already observed that $\bigl[ 2/(q^{1/2} - q^{-1/2})\bigr] \tau$
tends to the derivative operator $d/dx$ as $q \rightarrow 1$.  We therefore
expect the eigenfunctions of $\tau$ to provide a certain 
$q$-analog of the exponential.  Let us now determine these eigenfunctions.

Take the functions
$$
\psi_n(a, \cos \theta) = \bigl( aq^{(1-n)/2} e^{i\theta};q\bigr)_n\,
\bigl( aq^{(1-n)/2} e^{-i\theta};q\bigr)_n\ . \eqno(3.1)
$$
It is easy to check that they satisfy
$$
\tau\, \psi_n(a, \cos\theta) = a\,  q^{-n/2}(1-q^n)\,
\psi_{n-1}(a,\cos\theta)\ . \eqno(3.2)
$$
Following Ref.[14], define the function
$$
{\cal E}_q(x;a,b) = \sum_{n=0}^\infty {q^{n^2/4}\over(q;q)_n}\,
\psi_n(a,\cos\theta)\, b^n\ , \qquad x = \cos\theta\ . \eqno(3.3)
$$
It immediately follows from (3.2) that ${\cal E}_q(x;a,b)$ is an
eigenfunction of $\tau$:
$$
\tau \,{\cal E}_q(x;a,b) = abq^{-1/4}\, {\cal E}_q(x;a,b)\ . \eqno(3.4)
$$

The function ${\cal E}_q(x;a,b)$ is thus the 
$q$-analog of the exponential that we were looking for.  
Indeed as $q \rightarrow 1^-$, 
$\psi_n(a,\cos\theta) \rightarrow (1 + a^2 - 2ax)^n$ 
and $(1-q)^n/(q;q)_n \rightarrow 1/n!$.  Therefore,
$$
\lim_{q\rightarrow 1^-} {\cal E}_q \bigl( x;a,(1-q)b \bigr) = 
\exp \bigl[ (1 + a^2 - 2ax) b \bigr]\ , \eqno(3.5)
$$
and in particular, for $a=-i$,
$$
\lim_{q\rightarrow 1^-} {\cal E}_q\bigl( x; -i, (1-q)\,b/2 \bigr) = e^{ibx}\ .
\eqno(3.6)
$$
In the following, we shall need the value of ${\cal E}_q{(x;-i,b/2)}$
at $x=0$.  In this connection, one readily finds that
$$
i^n q^{n^2/4}\, \psi_n(-i,\cos\theta) \Big|_{\cos\theta = 0} =
\cases{(q;q^2)^2_{n/2}\ , &$n$ even\ ,\cr
0\ , &$n$ odd\ .}\eqno(3.7)
$$
As a result, one has that$^{14}$
$$
{\cal E}_q(0,-i,b/2) = \sum_{n=0}^\infty {(q;q^2)_n\over (q^2,q^2)_n}
\left(-{b^2\over4} \right)^n = 
{(-qb^2/4;q^2)_\infty\over (-b^2/4;q^2)_\infty}\ , \eqno(3.8)
$$
with the last equality obtained from the $q$-binomial theorem:
$$
\sum_{n=0}^\infty {(a;q)_n\over (q;q)_n}\, z^n = 
{(az;q)_\infty\over(z;q)_\infty}\ . \eqno(3.9)
$$

\vskip 2cm

{\bf 4. A GENERATING FUNCTION FOR THE CONTINUOUS}\hfill\break\indent
{\bf\phantom{4. }q-HERMITE POLYNOMIALS}
\bigskip

We shall now use the $q$-exponential we described in Section~3 as
well as the connection between the $q$-Heisenberg algebra and the
continuous $q$-Hermite polynomials to provide an algebraic derivation
of an identity obtained in Refs.[14,20]. Take the 
$q$-exponential ${\cal E}_q(x;-i,b/2)$ of the operator $x$ and act with
it on the representation space with basis $\{ f_n \}$
$$
{\cal E}_q \left(x; -i, b/2 \right)\, f_n = \sum_{k=0}^\infty U_{kn}(b)\,
f_k\ . \eqno(4.1)
$$
In the realization of Section~2, $f_n(x) = H_n(x \vert q)$ and so,
$f_0(x) = 1$.  Setting $U_{k0}(b) \equiv U_k(b)$, we thus have for
$n=0$:
$$
{\cal E}_q \left( x;-i, b/2 \right) = \sum_{k=0}^\infty U_k(b)\, 
H_k(x\vert q)\ . \eqno(4.2)
$$
To determine the matrix elements $U_k(b)$, act with $\tau$ on both
sides of (4.2).  Using (3.4) on the l.h.s. and (2.8) on the r.h.s., 
one finds
$$
-(ib/2)\,  q^{-1/4}\, {\cal E}_q \left(x; -i, b/2 \right) 
= \sum_{k=0}^\infty q^{k/2} (1-q^{-k})\,  U_k(b)\, H_{k-1}(x\vert q)\ . 
\eqno(4.3)
$$
Using again (4.2), the following two-term recurrence relation is
obtained
$$ 
(ib/2)\,  U_k(b) = q^{-(2k+1)/4} (1-q^{k+1})\, U_{k+1}(b)\ . \eqno(4.4)
$$
It has for solution
$$
U_k(b) = {q^{k^2/4}\over(q;q)_k} \left({ib\over2}\right)^k\, 
U_0(b)\ . \eqno(4.5)
$$
To determine $U_0(b)$, replace $U_k(b)$ in (4.3) by this last
expression and set $x=0$.  With the help of (3.8), one thus obtains
$$
{(-qb^2/4; q^2)_\infty\over(-b^2/4; q^2)_\infty} = U_0(b)\,
\sum_{k=0}^\infty {q^{{k}^2/4}\over(q;q)_k}
\left({ib\over 2}\right)^k\, H_k(0 \vert q)\ . \eqno(4.6)
$$
From the three-term recurrence relation (2.16) for $H_k(x \vert q)$
one easily gets:
$$
H_{2k}(0 \vert q)= (-1)^k (q;q^2)_k\ ,
\qquad H_{2k+1}(0 \vert q)=\, 0\ .\eqno(4.7)
$$
With this information, (4.6) is rewritten as
$$
{(-qb^2/4; q^2)_\infty\over(-b^2/4; q^2)_\infty} = U_0(b)\,
\sum_{k=0}^\infty {(q^2)^{k(k-1)/2}\over(q^2;q^2)_k}
\left({qb^2\over4}\right)^k\ . \eqno(4.8)
$$
Recalling (1.4), this yields
$$
U_0(b) ={1\over(-b^2/4;q^2)_\infty}\ . \eqno(4.9)
$$
Putting everything together, we get the following identity$^{14,20}$
$$
(-b^2/4;q^2)_\infty \ {\cal E}_q \left( x; -i,  b/2 \right) =
\sum_{k=0}^\infty {q^{{k}^2/4}\over(q;q)_k}
\left({ib\over2}\right)^k\, H_k(x \vert q )\ , \eqno(4.10)
$$
which provides a generating relation for the continuous $q$-Hermite
polynomials.

\vskip 2cm

{\bf 5. CONTINUOUS q-ULTRASPHERICAL POLYNOMIALS AND A}\hfill\break\indent
{\bf\phantom{5. }q-DEFORMATION OF E(3)}
\bigskip

The continuous $q$-Hermite polynomials $H_n(x \vert q)$ involve no
parameters.  We shall now work our way one step up in the $q$-Askey
scheme$^{2}$ and consider a class of continuous $q$-polynomials
depending on one parameter.  The continuous $q$-ultraspherical 
polynomials $C_n(x;q^m \vert q)$ have already been defined
in (1.18) through their relation with the Askey-Wilson polynomials. 
They can also be expressed as follows:$^{1,2}$ 
$$
C_n(x;q^m \vert q) =
	\sum_{k=0}^n {(q^m;q)_k\, (q^m; q)_{n-k}\over(q;q)_k\, (q;q)_{n-k}}\,
e^{i(n-2k)\theta}\ , \qquad x = \cos \theta\ , \eqno(5.1)
$$
from where it is directly seen that these polynomials are of definite
parity,
$$
C_n(-x; q^m \vert q) = (-1)^n C_n(x; q^m \vert q)\ . \eqno(5.2)
$$
We shall now be considering operators acting on functions of two
variables $x = (z + z^{-1})/2$ and $t$.

Let us first introduce the operators $J_+$, $J_-$ and $K$:
$$
\eqalignno{
&J_+={q^{1/2}\over1-q}\,  t \, T_t^{-1/2}\, \tau\ , &(5.3a)\cr
&J_-={q\over1-q}\, {1\over t}\,  T_t^{-1/2}\, \tilde \tau\ , &(5.3b)\cr
&K= q^{-1/2}\, T_t\ , &(5.3c)}
$$
where $\tau = (z - z^{-1})^{-1} (T_z^{1/2} - T_z^{-1/2})$ is again
the divided difference operator and$^{19}$
$$
\tilde \tau ={q^{-1/2}\over z - z^{-1}} \bigg[
{(1-z^2T_t)(1-qz^2T_t)\over z^2}\, T_z^{1/2} - z^2 
\bigg(1-{T_t\over z^2} \bigg) \bigg(1-q{T_t\over z^2} \bigg)\, T_z^{-1/2}\bigg]
\ . \eqno(5.4)
$$
If we set $t = e^{i\phi}$ and $K = q^{J_0}$, we note that in the
limit $q \rightarrow 1^-$:
$$
\eqalign{
&J_+\rightarrow {1\over2} {e^{i\phi}\over\sin\theta}\,
{\partial\over\partial\theta}\ ,\cr
&J_-\rightarrow -2 e^{-i\phi} \left( \sin \theta
{\partial\over\partial\theta} - 2i \cos \theta
{\partial\over\partial\phi} + \cos \theta \right)\ ,\cr
&J_0 \rightarrow -i \left({\partial\over\partial\phi} +{1\over2}\right)\ .}
\eqno(5.5)
$$
When $q \rightarrow 1^-$, we therefore have a realization of $sl(2)$ on
the 2-sphere.  It can be checked in fact, that the generators $J_+$,
$J_-$ and $K$ satisfy the relations
$$
[J_+, J_-] = {K-K^{-1}\over q^{1/2}-q^{-1/2}}\ ,\qquad 
K\,J_{\pm} =q^{\pm 1} J_\pm\, K\ , \eqno(5.6)
$$
and thus provide a realization of the quantum algebra 
${\cal U}_q\bigl(sl(2) \bigr)$.  The connection with the continuous
$q$-ultraspherical polynomials is now made by observing that in this
model, the basis vectors of the associated module are realized by
$$\eqalign{
Q_m^\ell (x,t) = {(q;q)_{\ell-m}\over(q^{2m};q)_{\ell - m}}\,
q^{m(\ell - m)/2}\, C_{\ell - m} & (x; q^m \vert q)\, t^m\ ,\cr  
&m \leq \ell\ ;\quad \ell ,\ m = 0,1,2 \dots\ .}\eqno(5.7)
$$
It is found that $J_\pm$ and $K$ act as follows on these basis
functions:
$$
\eqalignno{
&J_+ Q^\ell_m={q\over1-q}\,
{(1-q^{m-\ell})(1-q^{m+\ell})\over(1-q^{2m+1})(1-q^m)}\
Q^\ell_{m+1}\ &(5.8a)\cr
&J_- Q^\ell_m= -{q^{1-m}\over1-q}\, (1-q^{2m-1})(1+q^{m-1})\ Q_{m-1}^\ell\ ,
&(5.8b)\cr
&K\,Q^\ell_m= q^{m-1/2}\, Q^\ell_m\ . &(5.8c)}
$$
Actually, it turns out possible to realize an algebra larger than
${\cal U}_q\bigl(\sl(2)\bigr)$.  Indeed, the three independent
operators
$$
P_0 = x\ , \qquad P_+ = t\ , \qquad P_- ={1\over t}(1-x^2)\ , \eqno(5.9)
$$
are well-defined on the space spanned by the $Q_m^\ell$.  In the
limit $q \rightarrow 1^-$,  the six generators $J_+$, $J_-$, $J_0$, $P_+$,
$P_-$ and $P_0$ satisfy the commutation relations of the Euclidean
algebra in three dimensions,$^{21,22}$ which we shall
denote by $E(3)$.  It thus follows that this set of operators
defines a realization on the 2-sphere of a $q$-deformation of $E(3)$.
We shall now describe how $P_0$ and $P_+$ transform the basis
functions $Q^\ell_m$.

The action of $P_0$ is directly obtained from the three-term
recurrence relation of the continuous $q$-ultraspherical polynomials
$C_n(x; q^m \vert q)$:$^{1,2}$
$$
2x\,C_n(x; q^m \vert q) =
{1 - q^{n+1}\over1 - q^{m+n}}\, C_{n+1}(x;q^m \vert q) +
{1 - q^{2m + n - 1}\over1 - q^{m+n}}\, C_{n-1} (x;q^m \vert q)\ .\eqno(5.10)
$$
This yields
$$
P_0\, Q_m^\ell ={q^{-m/2}\over2} \left({1-q^{\ell+m}\over1-q^\ell}\right)\,
Q^{\ell + 1}_m +{q^{m/2}\over2} 
\left({1-q^{\ell- m}\over1-q^\ell}\right)\, Q^{\ell - 1}_m\ .\eqno(5.11)
$$
The action of $P_+ = t$ is also straightforwardly obtained.  In the
limit $q \rightarrow 1^-$, $P_+$ is a $sl(2)$ spin 1 operator and for
generic $q$, we know the $sl(2)$ and ${\cal U}_q\bigl(sl(2)\bigr)$ 
modules to be isomorphic.  Taking into account the fact that 
$KP_+ = q\,P_+ K$, we must have therefore
$$
P_+ Q^\ell_m = f_{\ell, m}\, Q^{\ell -1}_{m+1} + h_{\ell, m}\,
Q^{\ell +1}_{m+1}\ ,\eqno(5.12)
$$
with $f_{\ell, m}$ and $h_{\ell, m}$ coefficients to be determined. 
In terms of the polynomials $C_n(x; q^m \vert q)$, this means that
there must be an identity of the form
$$
\eqalign{
C_n(x; q^m \vert q) =& {q^{-m-1}\over(1-q^n)(1-q^{n-1})}\,
C_{n-2}(x;q^{m+1} \vert q)\, f_{n,m} \cr
&+{1\over(1 - q^{2m+n}) (1 - q^{2m+n+1})}\, C_n(x;q^{m+1} \vert q)\,h_{n,m}
\ .}\eqno(5.13)
$$
The constants $f_{n,m}$ and $h_{n,m}$ can now be obtained as follows.
The values of the $h_{n,m}$ are found first by equating the
coefficients of $x^n$ on both sides of (5.13) and using
$$
C_n (x; q^m \vert q) = 2^n\,{(q^m;q)_n\over (q;q)_n} \, x^n + \dots
\ .\eqno(5.14)
$$
With the result substituted back in (5.13), the $f_{n,m}$ are
determined by evaluating  at $x=0$.  The values of the continuous
$q$-ultraspherical polynomials at $x=0$ can be gotten for instance
from the recurrence relation (5.10) and are
$$
\eqalign{
&C_{2k}(0; q^m \vert q)= (-1)^k {(q^{2k}; q^2)_k\over(q^k; q^2)_k}\ ,\cr
&C_{2k+1}(0;q^m \vert q)=\,0\ .}
\eqno(5.15)
$$
If $n=2k$, setting $x=0$ in (5.13) will then immediately yield
$f_{2k, m}$.  In order to obtain $f_{2k+1,m}$, one first applies
$\tau$ to both sides of (5.13) knowing that
$$
\tau \, C_n(x;q^m \vert q) = -q^{-n/2} (1-q^m)\,  
C_{n-1}(x; q^{m+1}\vert q)\ , \eqno(5.16)
$$
before taking $x=0$.  One thus supplements (5.12) with
$$
\eqalignno{
&f_{\ell,m}= -q^{-(\ell - m)/2}\, 
{(1-q^m)(1-q^{\ell - m})(1-q^{\ell -m - 1}) q^{2m+1}\over
(1-q^{2m}) (1 - q^{2m+1})(1-q^\ell)}\ , &(5.17a)\cr
&h_{\ell,m}= q^{-(\ell - m)/2}\, 
{(1-q^m)(1-q^{\ell + m})(1-q^{\ell + m+1})\over
(1-q^{2m}) (1 - q^{2m+1})(1-q^\ell)}\ . &(5.17b)}
$$
The action of $P_-$ on the basis functions $Q^\ell_m$ is similarly
obtained and one has
$$
P_-\, Q^\ell_m = r_{\ell,m}\, Q^{\ell - 1}_{m-1} + s_{\ell, m}\,
Q_{m-1}^{\ell + 1}\ , \eqno(5.18)
$$
with 
$$
\eqalignno{
&r_{\ell , m}= q^{(\ell - m)/2}\, {(1-q^{2m-1})\over(1-q^{\ell+m-1})}
\left[ 1 + {(1-q^{2m-1})(1-q^{\ell -m+1})\over4(1-q^{m-1})(1-q^\ell)}\right]
\ , &(5.19a)\cr
&s_{\ell, m}= -q^{(\ell - 3m + 2)/2}\,
{(1-q^{2m-2})(1-q^{2m-1})\over4(1-q^{m-1})(1-q^\ell)}\ . &(5.19b)}
$$

\vskip 2cm

{\bf 6. A q-ANALOG OF THE FOURIER-GEGENBAUER EXPANSION}
\bigskip
 
We shall now make use of the results obtained in the previous
sections to derive a $q$-analog of the Fourier-Gegenbauer expansion
(1.21).  Take the ${\cal E}_q$ exponential of $P_0 = x$ and consider
its matrix elements in the basis $Q^\ell_m$:
$$
{\cal E}_q\left(x; -i, b/2\right)\, Q^\ell_m =
\sum_{\ell ' m'} U^{\ell', \ell}_{m',m}(b)\, Q^{\ell'}_{m'}\ .
\eqno(6.1)
$$
Recalling that $K\, Q^\ell_m = q^{m-1/2}\, Q^\ell_m$, we must have
$U^{\ell', \ell}_{m',m}(b) = U^{\ell', \ell}_{m,m}(b)\,  \delta_{mm'}$,
since $K = q^{-1/2}\, T_t$ and $P_0$ commute.  If we specialize (6.1)
to the case $m = \ell$, taking into account (5.11) and 
setting $U^{\ell+k, \ell}_{\ell,\ell}(b) \equiv W_\ell^k(b)$, we find
$$
{\cal E}_q \left(x; -i, b/2\right)\, Q^\ell_\ell = \sum^\infty_{k=0}
W_\ell^k(b)\, Q_\ell^{\ell+k}\ , \eqno(6.2)
$$
where from (5.7), we know that $Q_\ell^\ell = t^\ell$.  We shall see
that (6.2) contains the identity we are looking for.  We shall
exploit the representation of the $q$-deformation of $E(3)$ given in
Section~5, acting in turn with $P_+$ and $J_+$ on both
sides of (6.2), to derive the recursion relations that the matrix
elements of ${\cal E}_q(x; -i, b/2)$ obey and in the end, to determine
these $W_\ell^k(b)$.  Let us start with $P_+ = t$.  Since 
$t \, Q_\ell^\ell = Q_{\ell+1}^{\ell+1}$, we have on the one 
hand
$$
P_+\, {\cal E}_q \left( x; -i, b/2 \right)\, Q^\ell_\ell
={\cal E}_q \left( x; -i, b/2\right) \,Q^{\ell + 1}_{\ell + 1}
= \sum_{k=0}^\infty W_{\ell + 1}^k(b)\, Q^{\ell + k + 1}_{\ell +1}\ ,
\eqno(6.3)
$$
with the last equality following from (6.2).  On the other hand, 
$$
P_+\, {\cal E}_q \left(x; -i, b/2\right)\, Q_\ell^\ell = \sum_{k=0}^\infty
W_\ell^k(b)\, P_+ \, Q_\ell^{\ell + k}\ ; \eqno(6.4)
$$
using the action of $P_+$ on the basis functions $Q^\ell_m$, given in
(5.12), we can write
$$
P_+\, {\cal E}_q \left(x; -i, b/2 \right)\, Q_\ell^\ell = \sum_{k=0}^\infty
\bigl[ W_\ell^{k+2} (b) \,
f_{\ell +k+2,\ell} + W_\ell^k(b) \, h_{\ell + k,\ell} \bigr]\, 
Q_{\ell+1}^{\ell + k+1}\ , \eqno(6.5)
$$
with $f_{\ell,m}$ and $h_{\ell,m}$ as in (5.17).  Equations (6.3) and
(6.5) are then seen to imply the relation:
$$
W_{\ell +1}^k(b) = W_\ell^{k+2}(b)\, f_{\ell +k+2,\ell} 
+ W_\ell^k (b) \, h_{\ell+k,\ell}\ . \eqno(6.6)
$$

Now act similarly on both sides of (6.2) with $J_+$ to find:
$$
J_+\, {\cal E}_q(x; -i, b/2)\, Q_\ell^\ell = \sum_{k=0}^\infty
W_\ell^k(b) \, J_+ \, Q_\ell^{\ell +k}\ . \eqno(6.7)
$$
From the definition $(5.3a)$ of $J_+$ and the fundamental property
(3.4) of ${\cal E}_q$, it is immediate to check that
$$
\eqalign{
J_+\, {\cal E}_q\left(x; -i, b/2\right) \, Q_\ell^\ell
	&= -{i\over2}\,{q^{-(\ell - 1/2)/2}\over1-q}\, 
{\cal E}_q(x; -i,b/2) \, Q_{\ell +1}^{\ell +1}\cr
	&= -{i\over2}\, {q^{-(\ell - 1/2)/2}\over1-q}
\sum_{k=0}^\infty W_{\ell +1}^k(b) \, Q_{\ell +1}^{\ell +k+1}\ ,}
\eqno(6.8)
$$
with the second equality resulting from using (6.2) anew. The action
of $J_+$ on $Q_m^\ell$ is given in $(5.8a)$, and yields
$$
\sum_{k=0}^\infty W_\ell^k(b) \, J_+\, Q^{\ell + k}_\ell
= \sum_{k=0}^\infty W_\ell^{k+1}(b)
\left[{q(1-q^{-k-1})(1-q^{2\ell+k+1)}\over
(1-q)(1-q^{2\ell+1})(1+q^\ell)} \right] Q_{\ell +1}^{\ell + k+1}\ .
\eqno(6.9)
$$
These results imply:
$$
W_{\ell + 1}^k(b) ={2i\over b} q^{(2\ell+3)/4}\, 
{(1-q^{-k-1})(1-q^{2\ell + k + 1})\over(1-q^{2\ell+1})(1+q^\ell)}
\, W^{k+1}_\ell(b)\ . \eqno(6.10)
$$

It will now prove convenient to use instead of $W_\ell^k(b)$, the
functions $Y_\ell^k(b)$, which differs from $W_\ell^k(b)$  by a
multiplicative constant:
$$
W_\ell^k(b) = i^k q^{k(k-2\ell)/4}\,
{(q^{2\ell};q)_k\,(1-q^{k+\ell})\over(q;q)_k} \, Y_\ell^k(b)\ . \eqno(6.11)
$$
Upon replacing $W_{\ell+1}^k(b)$ in (6.6) by the r.h.s. of (6.10)
and inserting the explicit expressions for $f_{\ell + k + 2,\ell}$
and $h_{\ell+k,\ell}$, one arrives at the following relation
$$
{2\over b}\, (1-q^{k + \ell + 1}) \, Y_\ell^{k+1}(b)
	= q^{k+\ell+1} \, Y_\ell^{k+2} (b) + Y_\ell^k(b)\ . \eqno(6.12)
$$
It will determine $Y_\ell^k(b)$ up to a factor not indexed by $k$. 
Indeed, if we cast $Y_\ell^k(b)$ in the form
$$
Y_\ell^k(b) = {\cal Y}_{k+\ell}(b)\, V_\ell(b)\ , \eqno(6.13)
$$
the term $V_\ell(b)$ factors out from (6.12) which becomes
$$
{2\over b}\, (1-q^{s+1})\, {\cal Y}_{s+1}(b) 
	= q^{s+1}\, {\cal Y}_{s+2}(b) + {\cal Y}_s(b)\ , \qquad s = k+\ell\ .
\eqno(6.14)
$$
Comparing with (1.12), we recognize in (6.14), the three-term
recurrence relation of the $q$-Bessel function $J_s^{(2)}(b;q)$. 

As for the function $V_\ell(b)$, it can be obtained by substituting
(6.11) in (6.10) with
$$
Y_\ell^k(b) = J_{k+\ell}^{(2)} (b;q) \, V_\ell(b)\ . \eqno(6.15)
$$
This gives
$$
{b\over2(1-q^\ell)}\, V_{\ell + 1}(b) = V_\ell (b)\ ,
\eqno(6.16)
$$
which is readily found to imply that
$$
V_\ell(b) = \left({2\over b}\right)^{\ell} (q;q)_{\ell - 1}\, X(b)\ ,
\eqno(6.17)
$$
with $X(b)$, a function yet to be specified.  At this point, putting
everything together, equation (6.2) has been shown to entail the
expansion:
$$
\eqalign{
{\cal E}_q\left(x;i, b/2 \right) =
X(b)\, \sum_{k=0}^\infty i^k q^{{k^2}/4}(1-q^{k+\ell})\, &
{(q;q)_\infty\over(q^\ell;q)_\infty}\cr
&\times\left({2\over b}\right)^{\ell} J^{(2)}_{\ell+k}(b) \,
C_k(x;q^\ell \vert q)\ .}\eqno(6.18)
$$

The function $X(b)$ can now be determined by taking the limit 
$q^\ell \rightarrow 0$ in this equation.  
Note from (1.18) and (1.19), that the
continuous $q$-ultraspherical polynomials $C_n(x; q^\ell \vert q)$
reduce to the continuous $q$-Hermite polynomials when
$q^\ell \rightarrow 0$:
$$
C_n(q;0 \vert q) ={H_n(x \vert q)\over(q;q)_n}\ . \eqno(6.19)
$$
Recalling formula (1.13) for the large order behavior of the
functions $J_\nu^{(2)}(z;q)$, (6.18) is seen to yield
$$
{\cal E}_q\left(x; -i, b/2 \right) = X(b) \sum_{k=0}^\infty i^k
\left({b\over2}\right)^k\, {q^{k^2/4}\over(q;q)_K}\, H_k(x \vert q)\ ,
\eqno(6.20)
$$
when $q^\ell \rightarrow 0$.  Upon comparing (6.20), with the generating
relation (4.10) for the continuous $q$-Hermite polynomials which was
algebraically derived in Section~4, we see that
$$
V(b) = \left( -b^2/4; q \right)_\infty^{-1}\ . \eqno(6.21)
$$
We hence finally arrive at the following identity:
$$
\eqalign{
{\cal E}_q \left(x;-i, b/2 \right) =
&{(q;q)_\infty (2/b)^{\ell}\over(b^2/4;q^2)_\infty(q^\ell;q)_\infty}\cr
&\hskip 1cm\times\sum_{k=0}^\infty i^k q^{k^2/4}\,(1-q^{k+\ell})\, 
J_{\ell+k}^{(2)}(b;q) \, C_k(x;q^\ell \vert q)\ .}
\eqno(6.22)
$$
This expansion formula which was here constructively derived using
algebraic methods, was first obtained in Ref.[14].

\vskip 2cm

{\bf 7. CONCLUDING REMARKS}
\bigskip

We have undertaken in this paper a study of the continuous
$q$-orthogonal polynomials from the algebraic point of view.  Our
strategy has been to proceed in analogy with the Lie theory approach
to ordinary special functions and to construct models where the 
continuous $q$-polynomials occur as matrix elements or basis
vectors in representations of $q$-algebras.  Divided difference
operators have been seen to play a fundamental role in these models
as should have been expected considering the difference equations 
that the continuous $q$-polynomials obey.  It has also been
observed that the $q$-exponential ${\cal E}_q$ is the one to use in
this context to mimic the classical exponential map.

We examined first the continuous $q$-Hermite polynomials.  They have
no free parameters and thus lie at the bottom of the Askey tableau. 
They were shown to realize a representation space for the
$q$-Heisenberg algebra.  It was immediate to derive the generating 
relation (4.10) in this framework.  We then proceeded to
discuss similarly the continuous $q$-ultraspherical polynomials which
involve one free parameter.  These functions were shown to arise in
the basis vectors of a representation space of a $q$-deformation 
of $E(3)$.  Using these results and ``bootstrapping''
from the $q$-Hermite to the $q$-ultraspherical case, we derived
algebraically the $q$-analog of the Fourier-Gegenbauer expansion.

It would be of interest to interpret other $q$-continuous polynomials
in an analogous fashion and in particular, working our way up in the
Askey classification, to determine the symmetry algebras of the
higher families of continuous $q$-polynomials.  Work
along those lines is in progress.$^{23,24}$

\vskip 2cm

\centerline{\bf REFERENCES}
\bigskip

\item{1.}  Gasper, G. and Rahman, M., {\it Basic Hypergeometric Series},
(Cambridge University Press, Cambridge, 1990)
\smallskip
\item{2.} Koekoek, R. and Swarttouw, R.~F., The Askey-scheme of
hypergeometric orthogonal polynomials and its $q$-analogue, Report 94-05,
Delft University of Technology, 1994
\smallskip
\item{3.} Floreanini, R. and Vinet, L., $q$-Orthogonal polynomials
and the oscillator quantum group, Lett. Math. Phys. {\bf 22} (1991), 45-54
\smallskip
\item{4.} Floreanini, R. and Vinet, L., Quantum algebras and $q$-special
functions, Ann. Phys. {\bf 221} (1993) 53-70
\smallskip
\item{5.} Floreanini, R. and Vinet, L., Addition formulas for
$q$-Bessel functions, J. Math. Phys. {\bf 33}, 2984-2988 (1992)
\smallskip
\item{6.} Floreanini, R. and Vinet, L., Generalized $q$-Bessel functions,
Can. J. Phys. {\bf 72}, 345-354 (1994)
\smallskip
\item{7.} Floreanini, R. and Vinet, L., Using quantum algebras
in $q$-special function theory, Phys. Lett. A {\bf 170}, 21-28 (1992)
\smallskip
\item{8.} Floreanini, R. and Vinet, L., ${\cal U}_q(sl(2))$ and
$q$-special functions, Contemp. Math. {\bf 160}, 85-100 (1994)
\smallskip
\item{9.} Floreanini, R. and Vinet, L., An algebraic interpretation
of the $q$-hypergeometric functions, J. Group Theory Phys. {\bf 1}, 1-10
(1993)
\smallskip
\item{10.} Floreanini, R. and Vinet, L., Automorphisms of the $q$-oscillator
algebra and basic orthogonal polynomials, Phys. Lett. {\bf A180}, 
393-401 (1993)
\smallskip
\item{11.} Floreanini, R., Lapointe, L. and Vinet, L., A quantum algebra
approach to basic multivariable special functions, 
J. Phys. A {\bf 27}, 6781-6797 (1994)
\smallskip
\item{12.} Floreanini, R. and Vinet, L., $q$-Gamma and $q$-beta functions in 
quantum algebra representation theory, J. Comp. App. Math., to appear
\smallskip
\item{13.} Ismail, M. E. H., The basic bessel functions and polynomials,
SIAM J. Math. Anal. {\bf 12}, 454-468 (1981)
\smallskip
\item{14.} Ismail, M.~E.~H. and Zhang, R., Diagonalization of certain
integral operators, Adv. Math. {\bf 109}, 1-33 (1994)
\smallskip
\item{15.} Askey, R.~A. and Wilson, J.~A., 
Some basic hypergeometric polynomials that generalize Jacobi
polynomials, Memoirs Amer. Math. Soc., vol. 319, 1985
\smallskip
\item{16.} Rahman,~M., The linearization of the product of continuous 
$q$-Jacobi polynomials, Can. J. Math. {\bf 33}, 225-284 (1981)
\smallskip
\item{17.} Askey,~R. and Ismail, M.~E.~H., A generalization of 
the ultraspherical polynomials, Studies in Pure Mathematics, P. Erd\"os, ed.,
(Birkauser, Basel, 1983), pp. 55-78
\smallskip
\item{18.} {\it Higher Transcendental Functions}, 
Erdelyi,~A., Magnus,~W., Oberhettinger,~F., and Tricomi, F.~G., eds.,
(McGraw-Hill, New york, 1953)
\smallskip
\item{19.} Kalnins, E.~G. and Miller, W., Symmetry techniques for $q$-series:
Askey-Wilson polynomials, Rocky Mountain J. Math. {\bf 19}, 223-230 (1989)
\smallskip
\item{20.} Al-Salam,~W., A characterization of the Rogers 
$q$-Hermite polynomials, University of Alberta, preprint, 1994
\smallskip
\item{21.} Miller, W., {\it Symmetries and Separation of Variables},
(Addison-Wesley, Reading (Massachusetts), 1977)
\smallskip
\item{22.} Miller, W., Special functions and the complex 
Euclidean group in 3-space, I, J. Math. Phys. {\bf 9}, 1163-1175 (1968)
\smallskip
\item{23.} Floreanini, R., LeTourneux, J. and Vinet, L., An algebraic
interpretation of the continuous big $q$-Hermite polynomials,
Universit\'e de Montr\'eal, preprint, CRM-2246, 1995
\smallskip
\item{24.} Floreanini, R., LeTourneux, J. and Vinet, L., More on the
$q$-oscillator algebra and $q$-orthogonal polynomials,
Universit\'e de Montr\'eal, preprint, CRM-2248, 1995

\bye